%% file: Enumath2019.tex
\pdfoutput=1
\documentclass[epj]{webofc}
\usepackage[utf8]{inputenc}
\usepackage[varg]{txfonts}   
\usepackage{booktabs}
\usepackage{xcolor}
\definecolor{darkred}{rgb}{0.4,0.0,0.0}
\definecolor{darkgreen}{rgb}{0.0,0.4,0.0}
\definecolor{darkblue}{rgb}{0.0,0.0,0.4}
\usepackage[bookmarks,linktocpage,colorlinks,
    linkcolor = darkred,
    urlcolor  = darkblue,
    citecolor = darkgreen]{hyperref}
\usepackage{graphicx}

\usepackage{hyperref}
\usepackage{cleveref}
\usepackage{pgfplots}
  \usetikzlibrary{plotmarks}
  \usetikzlibrary{arrows.meta}

  \usepgfplotslibrary{patchplots}
  \usepackage{grffile}
\usepackage{pgfgantt}
\usepackage{pdflscape}
 \pgfplotsset{compat=newest}
 \pgfplotsset{plot coordinates/math parser=false}
 \pgfplotsset{grid style={dotted}}
 \pgfplotsset{scaled x ticks=false}
 \pgfplotsset{scaled y ticks=false}
 \newlength\figureheight
 \newlength\figurewidth

\input{myshortcuts.tex}

\graphicspath{{images/}{./images/tikz/}}
\wocname{}
\woctitle{}
\begin{document}
%
\selectlanguage{english}
\title{Model Order Reduction of Combustion Processes with Complex Front Dynamics}
\author{%
\firstname{Philipp}  \lastname{Krah}\inst{1,2} \and
\firstname{Mario} \lastname{Sroka}\inst{1} \and
\firstname{Julius} \lastname{Reiss}\inst{1}
}
\institute{%
Technische~Universit\"at~Berlin,
Institut~f\"ur~Str\"omungsmechanik~und~Technische~Akustik~(TUB),
M\"uller-Breslau-Straße~15,
10623~Berlin
\and
Institut~f\"ur~Mathematik, Technische~Universit\"at~Berlin,
Straße~des~17.~Juni~136,
10623~Berlin
}
\abstract{
In this work we present a data driven method, used to improve mode-based model order
reduction of transport fields with sharp fronts.
We assume that the original flow field $q(\vec{x},t)=f(\phi(\vec{x},t))$ can be reconstructed
by a front shape function $f$ and a level set function $\phi$.
The level set function is used to generate a local coordinate, which parametrizes
the distance to the front.
In this way, we are able to embed the local 1D description of the front
for complex 2D front dynamics with merging or splitting fronts, while seeking a low
rank description of $\phi$.
Here, the freedom of choosing $\phi$ far away from the front
can be used to find a low rank description of $\phi$ which accelerates
the convergence of $\norm{q- f(\phi_n)}$, when truncating $\phi$ after the $n$th mode.
We demonstrate the ability of this new ansatz for a 2D propagating flame with a moving
front.
}

\maketitle
\section{Introduction}
\label{sec:intro}

Nowadays combustion systems are studied by simulating the reactive Navier Stokes equations
with billions of degrees of freedom. The simulations are numerically expensive,
because computational resources scale with the number of degrees of freedom.
Therefore, model order reduction (MOR) techniques are desired to reduce
the number of relevant parameters, which describe the system.
Unfortunately, classical MOR techniques fail for these systems.
We aim for an improvement in this report.

Our method follows a data driven approach, where a set of $N$ snapshots,
$\{q(\x,t_i)\}_{i=1,\dots N}$, gathered during a
numerical simulation, is used to generate a reduced order model (ROM).
Here, most ROMs rely on separation of
variables:
    \begin{equation}
      q(\x,t) \approx \sum_{i=1}^n a_i(t) \psi_i(\x)
      \label{eq:linearcombiPOD}
    \end{equation}
in which the initial high fidelity field $q(\x,t)$ is represented by a set of
basis functions $\psi_i$ and their amplitudes $a_i$.
Based on \cref{eq:linearcombiPOD},
Petrov Galerkin and Galerkin methods (see for a review \cite{RowleyColoniusMurray2004})
project the original dynamics on a $n$-dimensional subspace spanned by the basis $\psi_i$.
However, the approximation error of the produced ROM crucially depends on the error made
in \cref{eq:linearcombiPOD}.
Unfortunately, in combustion systems transport dominated phenomena like moving
flame kernels with sharp gradients or traveling shock waves
critically slow down the convergence of \cref{eq:linearcombiPOD}.
For instance, this was numerically investigated by \cite{HuangDuraisamyMerkle2018} for reactive flows
and is theoretically quantified with help of the Kolmogorov $n$-width \cite{OhlbergerRave2015,GreifUrban2019}.

In order to handle transport dominated fields with sharp fronts,
we propose a nonlinear mapping of the solution manifold:
     \begin{equation}
       q(\x,t) = f(\phi)\qquad \text{s.~t. } \phi(\x,t) = \sum_{i=1}^n \tilde{a}_i(t) \tilde{\psi}_i(\x)\,.
       \label{eq:nonmap}
     \end{equation}
Here, the function $f:\mathbb{R}\to\mathbb{R}$ is simply the wave profile and $\phi$ describes the shift of the wave profile in time.

 This approach shares similar features as \cite{ReS2018,LeC2018,KaratzasBallarinRozza2019,RimMoeLeVeque2018}, but overcomes some of the problems associated with these methods.
Namely,
\cite{ReS2018,KaratzasBallarinRozza2019} depends on the choice of the shifts or their a priori knowledge
and becomes cumbersome if the topology of the moving object changes.
The latter is also the main drawback of \cite{RimMoeLeVeque2018}.
Furthermore, \cite{ReS2018,RimMoeLeVeque2018}
are mainly applied in one spatial dimension.
Although \cite{LeC2018} does not suffer from the aforementioned drawbacks, it lacks
physical interpretation and provides little insight of the underlying structure.\\

The report is organized as follows:
In the first two sections we introduce the basic idea for 1D and 2D advective systems
and explain the benefits of our concept.
One possible realization is provided in \cref{sec:FSR}, which is then applied to
a real life application in combustion -- the reduction of burning hydrogen with
complex front dynamics including topological changes (\cref{sec:2DburningFlame}).
Finally, we come back to \cite{ReS2018,LeC2018,Welper2017,RimMoeLeVeque2018} by comparing
the concept (\cref{sec:conclusion}).

\section{Basic Idea: 1D example - advective transport}
\begin{figure}[htp!]
  \centering
  \setlength\figureheight{0.35\linewidth}
  \setlength\figurewidth{0.45\linewidth}
  \includegraphics[width=0.5\textwidth]{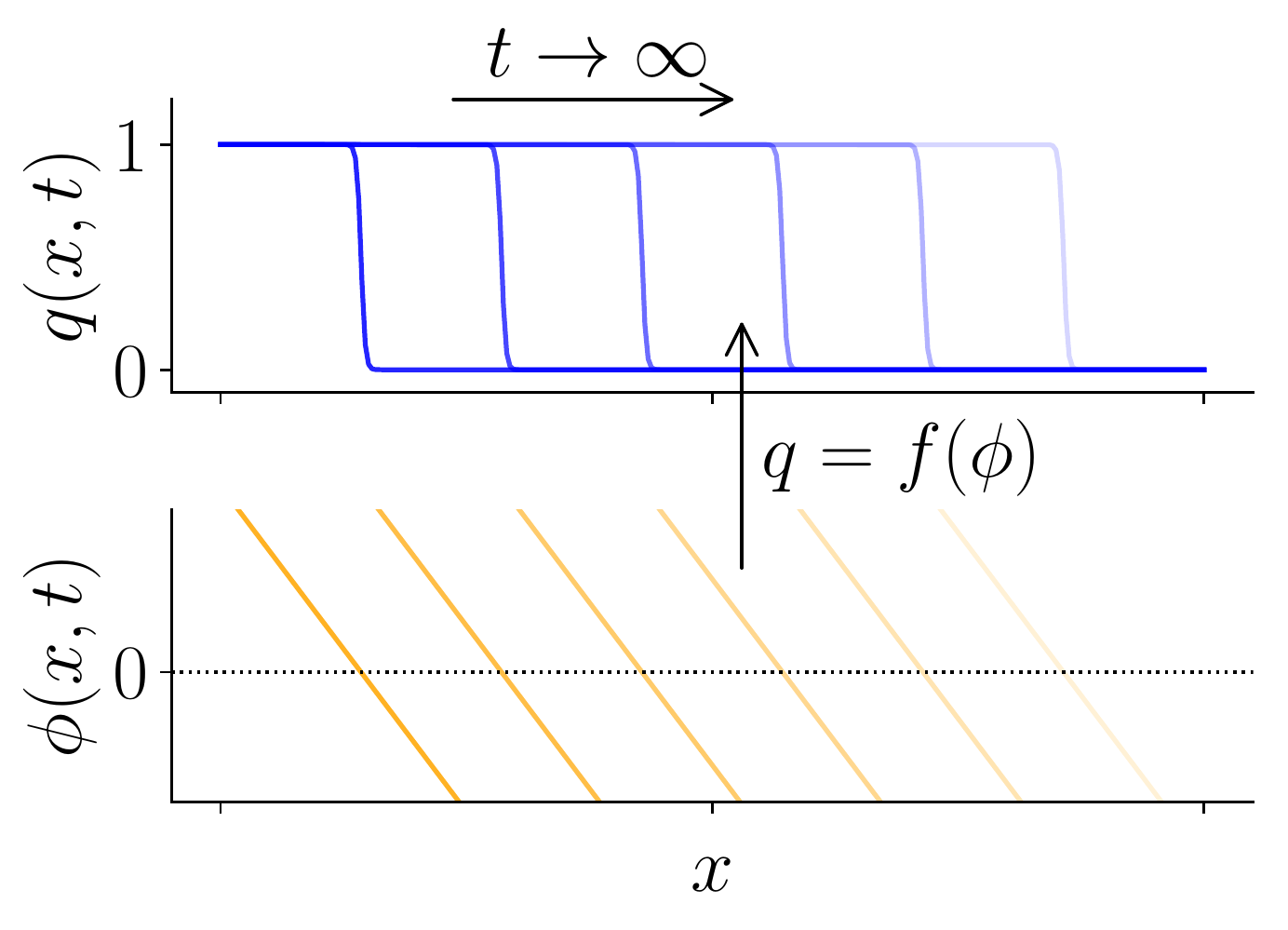}%
  \input{images/tikz/levelset1Dsigma.tex}
     \caption{Transported quantites $q$ and $\phi$ and singular values of the associated
       snapshotmatrix $X^q$ and $X^\phi$. Both functions share the same transport, since $q=f(\phi)$.
       However the transport of the sharp front is not well presented by a linear ansatz,
       \cref{eq:SVDX}, and therefore the singular values decay substantially slower as
     for the smooth field $\phi$. }
  \label{fig:1dlevelset}
\end{figure}

To motivate the proposed method we first consider a one dimensional problem.
A field $q(x,t)$ defined on $\mathcal{V}=[0,L]\times[0,T]$, with $L,T>0$ is given by
%
\begin{align}
  q(x,t)        &= f(\phi(x,t)),
\end{align}
where $f$ is a non-linear function and the auxiliary field
$  \phi(x,t)     = x-\Delta(t)$.
This describes an advective transport with trajectory path $\Delta(t)\colon[0,T]\to \mathbb{R}$,
in the most simple  example $\Delta (t) = c t$ with transport speed $c$.
The function $f$ is assumed to have a large gradient near $\phi = 0$.
In the examples we use
 \begin{equation}
   \label{eq:front}
  f_\lambda(\phi) =  (\tanh(\phi / \lambda) + 1)/2  \, ,
 \end{equation}
where $\lambda>0$ adjusts the front width.
Snapshots of the functions $q$ and $\phi$ are plotted for increasing time $t$ in
\cref{fig:1dlevelset}, left.
The corresponding snapshot matrices $X^\phi,\, X^q $ are defined as usual, $X^\alpha_{i,j} = \alpha(x_i,t_j) $.
A common approach to find a small representation of $X^\alpha$ is the truncated singular value decomposition (SVD)
\begin{equation}
  \label{eq:SVDX}
  {X^\alpha_n} = \sum_{k=1}^n \sigma_k \vec{u}_k \tran{\vec{v}_k} \, \widehat{=} \sum_{k=1}^n a_k(t) \psi_k(\vec{x})\,,
\end{equation}
which approximates $X^\alpha$ in the sense that the residuum $\Res=X^\alpha-X^\alpha_n$ is minimized.
We call the orthonormal basis $\{\psi_k(x_i)=(\vec{v}_k)_i\}_{k=1,\dots,n}$ spatial modes and
$\{a_k(t_i)=(\sigma_k\vec{u}_k)_i\}_{k=1,\dots,n}$ temporal coefficients.
As known from the Eckart Young Theorem \cite{EcY1936}, the approximation error $\norm{X^\alpha-X^\alpha_n}_2$
\footnote{{Note that for simplicity we will also use $\norm{\alpha-\tilde{\alpha}}_2$ for scalar functions $\alpha:\mathcal{V}\to \mathbb{R}$.
which we actually calculate as $\norm{X^\alpha-X^{\tilde{\alpha}}}_2$.}}
is given by the singular value $\sigma_{n+1}$, when truncating after the $n$th spatial mode.
The acceptable residuum is typically determined  by the target application.
A small number $n$ is desired in model order reduction as it governs the numerical cost of the reduced  model.

In \cref{fig:1dlevelset}, right we see the decay of the truncated singular value decomposition \cref{eq:SVDX}
of $X^q$ and $X^\phi$ is fundamentally different, even though $q$ is created from $\phi$ and both share the same advective
transport.
The failure of the SVD or POD to represent sharp transports is well known \cite{OhlbergerRave2015}.
In contrast to our example, $\phi$ can be represented by a linear combination of
two functions $\{x,1\}$.
Here, transport of the field $\phi$ is simply an amplitude change of the constant function.
This fact is not new and exploited by (\cite{ReS2018,Rim2018}).

With this ansatz we aim for a generalization of the method to higher spatial dimensions,
by  representing the movement of a front by an
auxiliary  field  $\phi$ which is of low rank
and a nonlinear mapping $f$ to recover the original field.
Thereby a locally one dimensional transport is implied.
However, we abstain from a global transport map between snapshots,
as this obstructs the application for topology changes.

\section{2D Example - Moving Disc}
\label{sec:2Dmovingdisc}

\begin{figure}[ht!]
	\centering
          \includegraphics[width=0.35\textwidth,trim=05 0 30 40,clip]{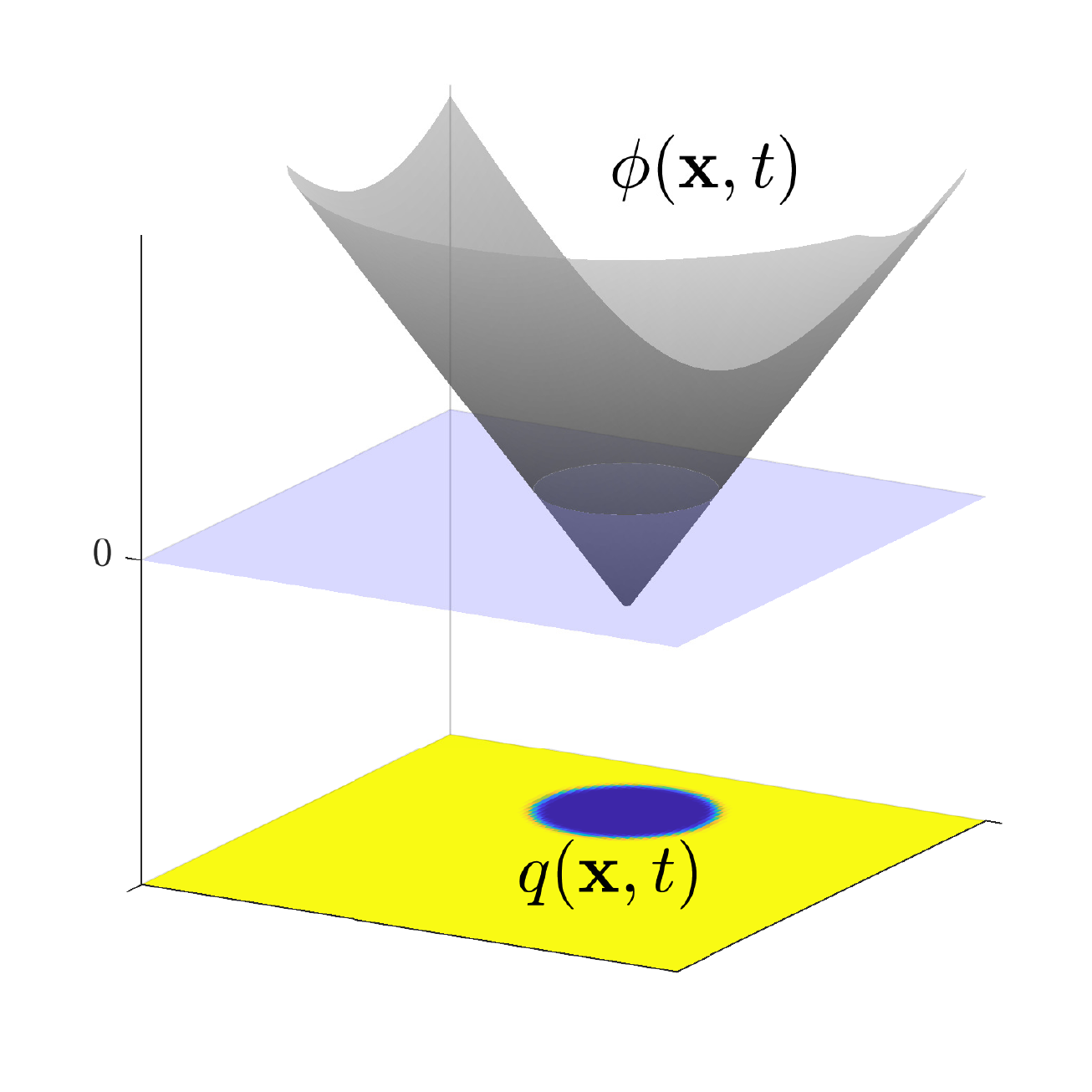}%
        \includegraphics[width=0.6\linewidth,trim=0 20 20 10,clip]{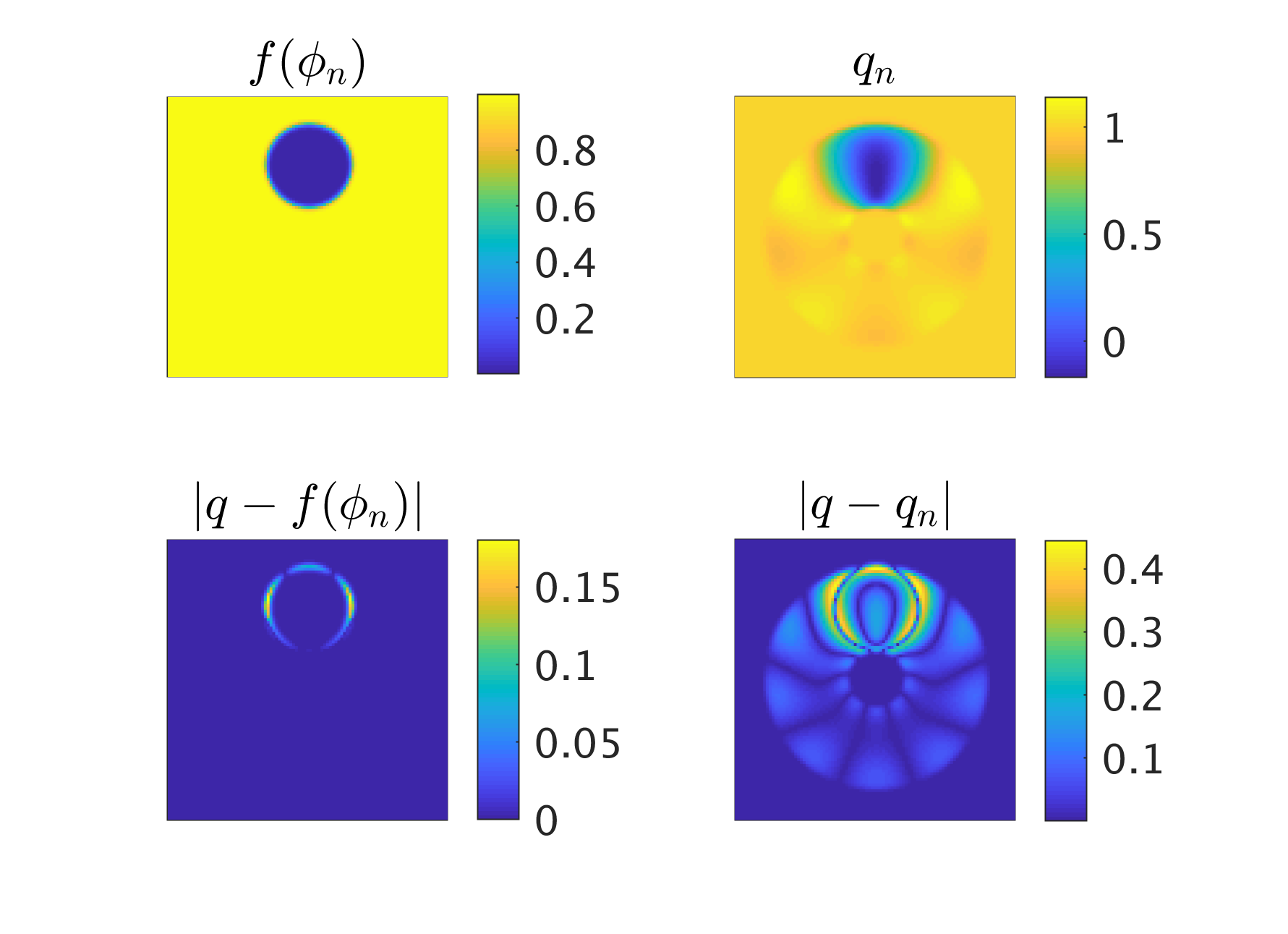}
	\caption{Left: visualization of the signed distance function. Right: Reconstruction of $q$ with $n=10$ modes.}
	\label{fig:translated-disc-reconst}
\end{figure}

The setting is now illustrated for a two dimensional problem
of a disc with radius $R=0.15L$, moving in a circle in a $[0,L]^2$ domain.
The translation of the disc is
parametrized by:
\begin{align}
  q(\vec{x},t)& = f(\phi(\vec{x},t))
        \quad\text{and}
\quad \phi(\vec{x},t)=\norm{\vec{x}-\vec{x}_0(t)}_2-R \label{eq:movDisc1}\\
	& \text{where } \vec{x}_0(t) =L
	\begin{pmatrix}
		0.5 + 1/4 \cos(2 \pi t)\\
		0.5 + 1/4 \sin(2\pi t)
	\end{pmatrix}\,,
	\label{eq:movDisc}
\end{align}
and $f$ is again the step function defined in \cref{eq:front}.
We sample 60 snapshots in a time interval $[0,1]$.
$\phi$ is the signed distance function shown in the left of \cref{fig:translated-disc-reconst},
which shall mimic the $\phi$ of the one dimensional example close to $\phi\approx0$.

\begin{figure}[htp!]
  \centering
  \setlength\figureheight{0.28\linewidth}
  \setlength\figurewidth{0.98\linewidth}
      \input{images/tikz/smoothed_phiv2.tex}
      \caption{Left:Graphical visualization of the smoothed signed distance function $\phi$
        defined in \cref{eq:movDisc1} and the paraboloid $\varphi(\vec{x},t)=\frac{1}{2R}(\norm{\vec{x}-\vec{x}_0(t)}_2^2-R^2)$ Left: A slice of
      $\phi$ and $\varphi$ at $t=0$ and $\vec{x}=(x,y_0),\, x\in [0,L]$.
     Comparison of the different truncation errors for the signed distance (center)
    and paraboloid (right).}
      \label{fig:smoothed_phi}
\end{figure}
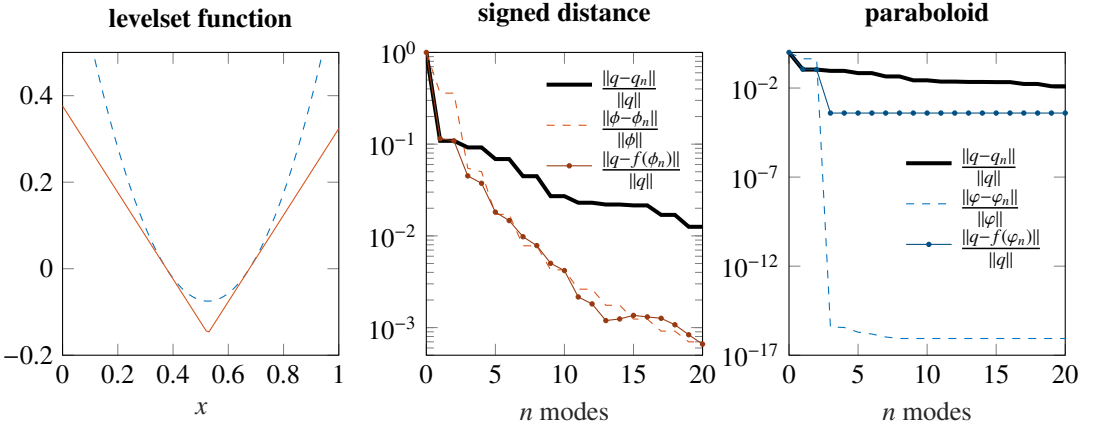

The original field $q$ is again reconstructed applying the SVD to
$\phi$ from which the approximation $\tilde q = f(\phi_n)$ is obtained.
\Cref{fig:translated-disc-reconst} shows the comparison between
the reconstruction using $f(\phi_n)$ and the naive POD approach using $q_n$ for snapshot $t=1/4$ with $n=10$ modes.
The results show not only a reduction in the overall error but also that the basic
structure of the moving disc is recovered.
The latter is already the case for a small number of modes.
While the example shows that the concept works well for 2D problems,
we show that the choice of $\phi$ has a strong effect on the possible
reduction.
By replacing $\phi$ with a paraboloid $\varphi(\vec{x},t)=\frac{1}{2R}(\norm{\vec{x}-\vec{x}_0(t)}_2^2-R^2)$
  in \cref{eq:movDisc1}, it is possible to obtain a much smaller error with less spatial modes.
  A one dimensional profile of the two functions $\phi,\varphi$ is plotted in the left of \cref{fig:smoothed_phi}.
 Note that $\phi$ and $\varphi$ have the same zero-level and their gradients are identical at the zero-level for all times $t$.
 Therefore $f(\varphi)\approx q=f(\phi)$ is a good approximation for small widths $\lambda$.
 In contrast to $\phi$, the total error
 $\norm{q-f(\varphi_n)}_2$ is reduced to its minimum with only three modes, because $\varphi$ can be
 represented by $\{x^2 + y^2 +R^2,x,y\}$.
  This is shown in the second and third column of \cref{fig:smoothed_phi}.
  From the ansatz $q\approx f(\varphi)$ we can deduce two different errors contributing to the total error:
\begin{equation}
  \label{eq:error}
  \norm{q-f(\varphi_n)} \le \underbrace{\norm{q-f(\varphi)}}_{=\Delta f}+\norm{ f'(\varphi)\Res} +\ord{\norm{\Res^2}}.
\end{equation}
The truncation error of the SVD is $\Res=\varphi-\varphi_n$ and the approximation error of the data is $\Delta f$.
Consequently, for vanishing approximation error, the truncation error
$\norm{f'(\varphi)\Res}\le\norm{f'(\varphi)}\sigma_{n+1}$ bounds the total error in the two norm.
Therefore, the total error scales with the decaying singular values of $\phi$.
This behaviour can be seen in \cref{fig:smoothed_phi}. While for $\phi$ the relative
truncation error aligns with the total error (middle), the error of $\varphi$ in
the right plot is dominated by the approximation error $\Delta f\approx 10^{-3}$.

From this example, we see that the ansatz is a good candidate for a low rank optimization of $\varphi$.
In contrast to areas where $f'(\varphi)\ne0$ and the field $\varphi$ has to mimic a signed distance to the
front, it can be chosen to minimize the truncation error far away from the zero level where $f'(\phi)=0$.
Additionally, in an optimization procedure one could relax the assumption of constant
front width by imposing appropriate conditions on the slope of $\varphi$ close to the zero level.
However this does not lie within the scope of this work.

 \section{Front Transport Reconstruction (FTR)}
\label{sec:FSR}
Now, we proceed to extent the idea to numerical data.
Here, only the field $q$ is known and the auxiliary field $\phi$ and
the front shape function $f$ need to be determined.
 For this, we assume that the front location can be calculated using threshold search
 of the relevant variables.

As proof of our concept, we compute $\phi$ as a two dimensional signed distance function,
because it is easy to compute and can be directly interpreted as a local 1D coordinate system.
This is a special choice for $\phi$ which is likely to be sub-optimal as was shown in the previous section. 
The zero-level curve $\segCurve$ of $\phi$ is determined by a threshold search with threshold $q_\segCurve$.
The discrete contour line $\segCurve$ was sampled at points where $q$ had the  value $q_\segCurve$ 
on any vertical or horizontal gridline of our computational
mesh. A linear interpolation of $q$ between the grid points is used to determine the crossing. 
The distance $ d_{\segCurve}(\x)$ is calculated as the minimal distance to all  sections of this curve, assumed to be linear between two points.
%
The sign of $\phi$ is negative if $ q(\x)< q_{\segCurve}$, and positive otherwise.
With the described procedure we determine the signed distance function $\phi(\vec{x},t_i)$
for every $t_i=i\Delta t$\,, $i=1,\dots,N$.

At this point a value of $\phi$ and $q$ is available at every grid point from which the front shape function
$f$ is to be determined such that $ q= f(\phi)$.
This is complicated by the fact that such relation is approximate and only discrete
values are available.
From the computed signed distance function we choose all grid points
$\hat{\phi}_l=\phi(x_{i_l},y_{j_l},t_{i_l})$ with  $\abs{\hat{\phi}_l}\le \Delta \phi$
on vertical, horizontal or diagonal lines which cross $\segCurve$
as support of the samples $\hat{q}_{l}=q(x_{i_l},y_{j_l},t_{k_l})$.
The sample vectors $(\hat{\phi}_l,\hat{q}_l)$ are then interpolated on a predefined support set $\phi_1,\dots,\phi_M$
which is used to find the corresponding interpolated values $f_1,\dots,f_M$ minimizing
the difference between $\hat{q}_l$ and $f(\hat{\phi}_l)$.

\section{2D Example - Application to Combustion}
\label{sec:2DburningFlame}

In this section we show that the described procedure is capable of reconstructing flow dynamics with inherent
two dimensional transport including  changing typology, which is difficult  for methods building on a
mapping between snapshots to remove transports.

\begin{figure}
  \setlength\figureheight{0.28\linewidth}
  \setlength\figurewidth{0.3\linewidth}
      \includegraphics[width=0.6\linewidth]{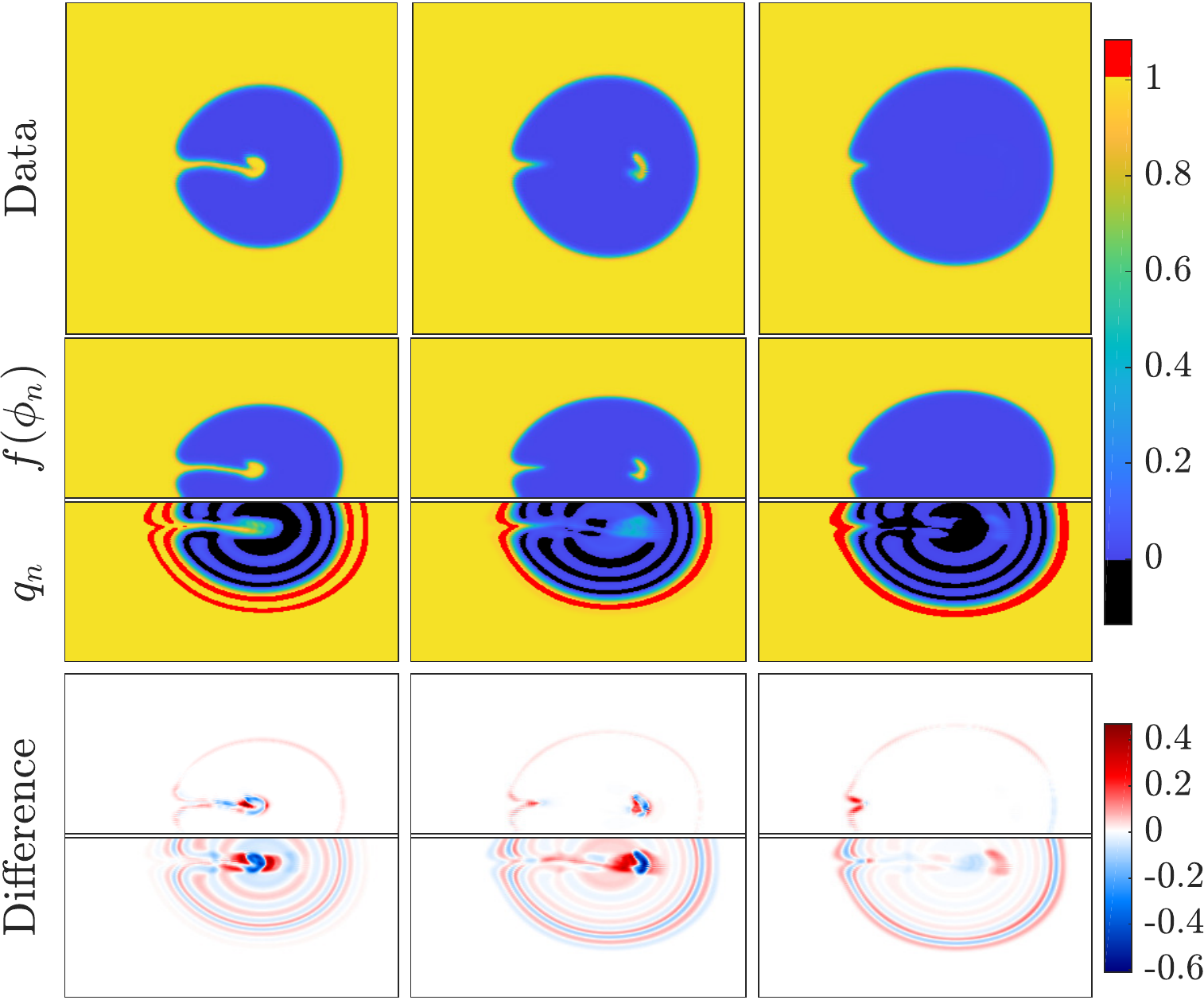}%
      \input{images/tikz/L2error1.tex}
      \caption{Comparison of $q_n$ (POD) and $f(\phi_n)$ (levelset method).
                 Left: Direct comparison of the snapshots $t/\Delta t=  21, 26, 31$
                 with $q$ plotted in the first row, the approximations $f(\phi_n)$
                 , $q_n$ in the second row and the difference between the data and its
                 approximation in the last row.
                  Please note that the images in the lower rows contain only the
                  fractions of the full snapshot that are relevant for our comparison.
              Right: Relative error in the two norm. }
      \label{fig:reconstr_flame}
\end{figure}

The configuration of a flame kernel interacting with a vortex pair mimics turbulence flame interaction.
Our data set consists of 40 snapshots derived from a 2D simulation of the reactive Navier Stokes equations.
For our purpose we restrict the reconstruction on the normalized mass fraction of hydrogen $Y_{\mathrm{H}_2}$.
 The simulation was tuned such that a
vortex pair moves towards burning H$_2$ and mixes unburned ($Y_{\mathrm{H}_2}=1$) with burned gas ($Y_{\mathrm{H}_2}=0$),
such that a small bubble of unburned gas detaches into the burned area.
The time evolution is visualized for some selected snapshots in the left of \cref{fig:reconstr_flame}.
As seen from \cref{fig:reconstr_flame}, the $Y_{\mathrm{H}_2}$ snapshots contain a very interesting
structure, in which the front changes along its contour line and even the topology of the
line changes -- splitting from one curve at $t/\Delta t = 24$ into two curves $t/\Delta t=29$ and
then back to a single curve at $t/\Delta t = 34$.

Applying the described procedure with a threshold of $q_\segCurve=0.14$, we achieve
promising results when comparing our method with a POD approximation in
\cref{fig:reconstr_flame} using 10 modes.
For this specific data the threshold $q_\segCurve$ was chosen to be rather small,
in order to resolve the tail of the incoming bubble.
The overall relative approximation error is decreased by a factor of three.
More important, our approximation preserves the physical structure of the data.
The POD does not respect the allowed physical range of $0\le Y_{\mathrm{H}_2}\le 1$ and
shows staircasing, i.e. replaces a front by several fractional fronts.
No sensible physical description can be expected from this structure.
The new method, in contrast, has well defined fronts and respects the physical range,
since $f$ is by construction restricted to the range of the input data $q$.
\section{Discussion and Conclusion}
\label{sec:conclusion}

We presented a concept for modal decomposition of transported fronts.
It builds on representing the original field $q$ by an auxiliary field $\phi$ and a non-linear function $f$ in such a way that the
$\phi$ has a better low rank description than $q$.
In a numerical example, $\phi$ was taken as a signed distance function with
the front as zero level and $f$ describing the front shape.
It is evident, that this choice is in general not optimal,
since moving kinks in $\phi$ yield a slow decay of singular values.

This approach can be interpreted as embedding a local one dimensional
coordinate into a multidimensional domain, orthogonal to the front.
A transport in this direction is simply an additive term for $\phi$.
The time dependent shift to compensate a transport in one dimension is a special case of this approach.
This induces a transport map, similar to \cite{Rim2018,ReS2018},
with the important difference that there is a local but no global
one to one mapping, by which topology changes are permitted.
A different perspective is the comparison with neural networks.
This linear combination to construct a low rank representation $\phi$
with the application of a non linear function can be
seen as a one layer network with a special activation function $f$.
A recent work uses level sets to handle geometry changes in which shares some technical
aspects with the current work \cite{KaratzasBallarinRozza2019}.

For full practical applicability, improvements are needed but near at hand.
To improve the approximation error $f$ and $\phi$ should be minimized based on \cref{eq:error}
and a more general ansatz should be used to allow a changing front shape.

\section*{Acknowledgments}
\label{Acknowledgments}
This work is supported by the Deutsche Forschungsgemeinschaft (DFG, German Research
Foundation) - 384950143/ GRK2433 and Project 200291049 - SFB 1029

\bibliography{refs.bib}

\end{document}

%% file: myshortcuts.tex
\usepackage{amsmath,amssymb,mathtools}

\renewcommand{\vec}[1]{\mathbf{#1}}

\newcommand{\tran}[1]{#1^{\intercal}}
\newcommand{\abs}[1]{\left | #1 \right |}
\newcommand{\ord}[1]{\mathcal{O}\left( #1 \right) }

\renewcommand{\vec}[1]{\mathbf{#1}}
\newcommand{\norm}[1]{\left\Vert #1 \right\Vert}
\newcommand{\Res}{\mathcal{R}}

\newcommand{\levelCurve}{\mathcal{C}}
\newcommand{\segCurve}{{\levelCurve_0}}

\newcommand{\x}{\vec{x}}

%% file: images/tikz/levelset1Dsigma.tex
\begin{tikzpicture}

\definecolor{color0}{rgb}{1,0.647058823529412,0}

\begin{axis}[
height=\figureheight,
legend cell align={left},
legend style={at={(0.5,0.6)}, anchor=north, draw=none},
log basis y={10},
tick align=outside,
tick pos=left,
title={${\sigma_n(X)}/{\Vert X\Vert_2}$},
width=\figurewidth,
x grid style={white!69.01960784313725!black},
xlabel={Mode \(\displaystyle n\)},
xmajorgrids,
xmin=0, xmax=20,
xminorgrids,
xtick style={color=black},
y grid style={white!69.01960784313725!black},
ymajorgrids,
ymin=4.54213436950621e-18, ymax=6.69645204306777,
yminorgrids,
ymode=log,
ytick style={color=black}
]
\addplot [blue, mark=*, mark size=1, mark options={solid}]
table {%
0 1
1 0.333307703774789
2 0.199954027121491
3 0.142791811095647
4 0.111027012985494
5 0.0908068272042107
6 0.0768032340419105
7 0.0665299577867899
8 0.0586707576941898
9 0.0524636522024691
10 0.0474369849968113
11 0.0432831983625682
12 0.0397931956823475
13 0.0368198763784882
14 0.0342567600498823
15 0.0320248853257144
16 0.0300644729534269
17 0.0283294474291522
18 0.0267837355609547
19 0.0253987040912931
20 0.0241513474298762
21 0.0230229812714542
22 0.0219982847067642
23 0.0210645870127596
24 0.020211329204153
25 0.0194296523606727
26 0.0187120792278299
27 0.0180522653322048
28 0.0174448025186014
29 0.0168850624490217
30 0.0163690708692308
31 0.0158934057818367
32 0.0154551143525037
33 0.0150516446111154
34 0.0146807889238899
35 0.014340636896266
36 0.0140295358836925
37 0.013746057681911
38 0.0134889702748358
39 0.0132572135707566
40 0.0130498795780921
41 0.0128661931241879
42 0.0127054977392726
43 0.0125672391273757
44 0.0124514240610197
45 0.0123563603073238
46 0.0122830307689976
47 0.0122308092701341
48 0.0121995376914025
49 0.00557113083247043
};
\addlegendentry{$X^q_{ij}=q(x_i,t_j)$}
\addplot [color0, dashed, mark=triangle*, mark size=1, mark options={solid,rotate=90}]
table {%
0 1
1 0.983468842746643
2 2.98409190362933e-16
3 2.39087578108502e-16
4 1.71640175526394e-16
5 1.06455300209217e-16
6 9.95076345691556e-17
7 9.95076345691556e-17
8 9.95076345691556e-17
9 9.95076345691556e-17
10 9.95076345691556e-17
11 9.95076345691556e-17
12 9.95076345691556e-17
13 9.95076345691556e-17
14 9.95076345691556e-17
15 9.95076345691556e-17
16 9.95076345691556e-17
17 9.95076345691556e-17
18 9.95076345691556e-17
19 9.95076345691556e-17
20 9.95076345691556e-17
21 9.95076345691556e-17
22 9.95076345691556e-17
23 9.95076345691556e-17
24 9.95076345691556e-17
25 9.95076345691556e-17
26 9.95076345691556e-17
27 9.95076345691556e-17
28 9.95076345691556e-17
29 9.95076345691556e-17
30 9.95076345691556e-17
31 9.95076345691556e-17
32 9.95076345691556e-17
33 9.95076345691556e-17
34 9.95076345691556e-17
35 9.95076345691556e-17
36 9.95076345691556e-17
37 9.95076345691556e-17
38 9.95076345691556e-17
39 9.95076345691556e-17
40 9.95076345691556e-17
41 9.95076345691556e-17
42 9.95076345691556e-17
43 9.95076345691556e-17
44 9.95076345691556e-17
45 9.95076345691556e-17
46 9.95076345691556e-17
47 9.95076345691556e-17
48 9.95076345691556e-17
49 3.0416184978568e-17
};
\addlegendentry{$X^\phi_{ij}=\phi(x_i,t_j)$}
\end{axis}

\end{tikzpicture}

%% file: images/tikz/smoothed_phiv2.tex
%
\definecolor{mycolor1}{rgb}{0.00000,0.44700,0.74100}%
\definecolor{mycolor2}{rgb}{0.85000,0.32500,0.09800}%
\begin{tikzpicture}

\begin{axis}[%
width=0.262\figurewidth,
height=\figureheight,
at={(0.689\figurewidth,0\figureheight)},
scale only axis,
xmin=0,
xmax=20,
xlabel style={font=\color{white!15!black}},
xlabel={$n$ modes},
ymode=log,
ymin=1e-17,
ymax=1,
yminorticks=true,
axis background/.style={fill=white},
title style={font=\bfseries},
title={paraboloid},
legend style={at={(0.97,0.5)}, anchor=east, legend cell align=left, align=left, fill=none, draw=none},
scaled ticks=false, xticklabel style={/pgf/number format/fixed},yticklabel style={/pgf/number format/fixed}
]
\addplot [color=black, line width=1.5pt]
  table[row sep=crcr]{%
0	1\\
1	0.108715760217601\\
2	0.108715634236489\\
3	0.0920641111130334\\
4	0.0920640847955109\\
5	0.068791913048795\\
6	0.0687918193595006\\
7	0.0447255749252609\\
8	0.044725562332172\\
9	0.0270961356129575\\
10	0.027096041422363\\
11	0.0229775396649234\\
12	0.0229774237831645\\
13	0.0219687207092515\\
14	0.0219685229489843\\
15	0.0214905919511001\\
16	0.0214903863358998\\
17	0.0169406684748452\\
18	0.0169404955096078\\
19	0.0125606794962384\\
20	0.0125602504887935\\
21	0.012095808160281\\
22	0.012095753375457\\
23	0.0114919120191121\\
24	0.0114915874902618\\
25	0.0114791822764769\\
26	0.0114791164300901\\
27	0.00942848486431508\\
28	0.00942818463158277\\
29	0.00753288540033795\\
30	0.00753278297038982\\
31	0.00752777202203293\\
32	0.0075272738516459\\
33	0.00720851011526472\\
34	0.00720817840874839\\
35	0.00712811894643316\\
36	0.00712798767233169\\
37	0.00598874820673185\\
38	0.00598869108086194\\
39	0.0051289427089822\\
40	0.00512891965034824\\
41	0.00503424481958879\\
42	0.00503424452216082\\
43	0.00494515674123981\\
44	0.00494492780082036\\
45	0.00487137514777922\\
46	0.00487037432500624\\
47	0.00421444940822556\\
48	0.00421409967359719\\
49	0.00395873561982647\\
50	0.00395846526079949\\
51	0.00386362980513069\\
52	0.00386353684236286\\
53	0.00377082854130412\\
54	0.00377030557373398\\
55	0.00372967042540982\\
56	0.00372896245992916\\
57	0.00351333695898237\\
58	0.00351271540544044\\
59	0.00330378411420591\\
};
\addlegendentry{$\frac{\Vert q - q_n\Vert}{\Vert q \Vert}$}

\addplot [color=mycolor1, dashed]
  table[row sep=crcr]{%
0	1\\
1	0.436930895045629\\
2	0.436887632419151\\
3	3.8193175378319e-16\\
4	3.61094984130645e-16\\
5	1.96856041779954e-16\\
6	1.54143069352301e-16\\
7	1.05611075450375e-16\\
8	8.67295662566665e-17\\
9	8.67295662566665e-17\\
10	8.67295662566665e-17\\
11	8.67295662566665e-17\\
12	8.67295662566665e-17\\
13	8.67295662566665e-17\\
14	8.67295662566665e-17\\
15	8.67295662566665e-17\\
16	8.67295662566665e-17\\
17	8.67295662566665e-17\\
18	8.67295662566665e-17\\
19	8.67295662566665e-17\\
20	8.67295662566665e-17\\
21	8.67295662566665e-17\\
22	8.67295662566665e-17\\
23	8.67295662566665e-17\\
24	8.67295662566665e-17\\
25	8.67295662566665e-17\\
26	8.67295662566665e-17\\
27	8.67295662566665e-17\\
28	8.67295662566665e-17\\
29	8.67295662566665e-17\\
30	8.67295662566665e-17\\
31	8.67295662566665e-17\\
32	8.67295662566665e-17\\
33	8.67295662566665e-17\\
34	8.67295662566665e-17\\
35	8.67295662566665e-17\\
36	8.67295662566665e-17\\
37	8.67295662566665e-17\\
38	8.67295662566665e-17\\
39	8.67295662566665e-17\\
40	8.67295662566665e-17\\
41	8.67295662566665e-17\\
42	8.67295662566665e-17\\
43	8.67295662566665e-17\\
44	8.67295662566665e-17\\
45	8.67295662566665e-17\\
46	8.67295662566665e-17\\
47	8.67295662566665e-17\\
48	8.67295662566665e-17\\
49	8.67295662566665e-17\\
50	8.67295662566665e-17\\
51	8.67295662566665e-17\\
52	8.67295662566665e-17\\
53	8.67295662566665e-17\\
54	8.67295662566665e-17\\
55	8.67295662566665e-17\\
56	8.67295662566665e-17\\
57	8.67295662566665e-17\\
58	5.12180668652929e-17\\
59	1.62766191257076e-17\\
};
\addlegendentry{$ \frac{\Vert\varphi-\varphi_n\Vert}{\Vert \varphi \Vert}$}

\addplot [color=black!30!mycolor1, mark size=0.8pt, mark=*, mark options={solid, fill=mycolor1, black!30!mycolor1}]
  table[row sep=crcr]{%
0	1\\
1	0.114746569351737\\
2	0.108715634236489\\
3	0.000394256743858815\\
4	0.000394256743858815\\
5	0.000394256743858814\\
6	0.000394256743858814\\
7	0.000394256743858813\\
8	0.000394256743858813\\
9	0.000394256743858811\\
10	0.000394256743858811\\
11	0.000394256743858808\\
12	0.000394256743858809\\
13	0.000394256743858809\\
14	0.000394256743858807\\
15	0.000394256743858802\\
16	0.000394256743858802\\
17	0.0003942567438588\\
18	0.000394256743858801\\
19	0.000394256743858801\\
20	0.0003942567438588\\
21	0.000394256743858799\\
22	0.000394256743858797\\
23	0.000394256743858796\\
24	0.000394256743858795\\
25	0.000394256743858796\\
26	0.000394256743858795\\
27	0.000394256743858796\\
28	0.000394256743858796\\
29	0.000394256743858793\\
30	0.000394256743858792\\
31	0.000394256743858791\\
32	0.000394256743858791\\
33	0.000394256743858791\\
34	0.000394256743858792\\
35	0.000394256743858792\\
36	0.000394256743858792\\
37	0.000394256743858792\\
38	0.000394256743858794\\
39	0.000394256743858794\\
40	0.000394256743858793\\
41	0.000394256743858794\\
42	0.000394256743858794\\
43	0.000394256743858794\\
44	0.000394256743858793\\
45	0.000394256743858794\\
46	0.000394256743858794\\
47	0.000394256743858794\\
48	0.000394256743858794\\
49	0.000394256743858794\\
50	0.000394256743858793\\
51	0.000394256743858793\\
52	0.000394256743858793\\
53	0.000394256743858794\\
54	0.000394256743858794\\
55	0.000394256743858795\\
56	0.000394256743858796\\
57	0.000394256743858796\\
58	0.000394256743858797\\
59	0.000394256743858797\\
};
\addlegendentry{$\frac{\Vert q - f(\varphi_n)\Vert}{\Vert q \Vert}$}

\end{axis}

\begin{axis}[%
width=0.262\figurewidth,
height=\figureheight,
at={(0.345\figurewidth,0\figureheight)},
scale only axis,
xmin=0,
xmax=20,
xlabel style={font=\color{white!15!black}},
xlabel={$n$ modes},
ymode=log,
ymin=0.0005,
ymax=1,
yminorticks=true,
axis background/.style={fill=white},
title style={font=\bfseries},
title={signed distance},
legend style={legend cell align=left, align=left, fill=none, draw=none},
scaled ticks=false, xticklabel style={/pgf/number format/fixed},yticklabel style={/pgf/number format/fixed}
]
\addplot [color=black, line width=1.5pt]
  table[row sep=crcr]{%
0	1\\
1	0.108715760217601\\
2	0.108715634236489\\
3	0.0920641111130334\\
4	0.0920640847955109\\
5	0.068791913048795\\
6	0.0687918193595006\\
7	0.0447255749252609\\
8	0.044725562332172\\
9	0.0270961356129575\\
10	0.027096041422363\\
11	0.0229775396649234\\
12	0.0229774237831645\\
13	0.0219687207092515\\
14	0.0219685229489843\\
15	0.0214905919511001\\
16	0.0214903863358998\\
17	0.0169406684748452\\
18	0.0169404955096078\\
19	0.0125606794962384\\
20	0.0125602504887935\\
21	0.012095808160281\\
22	0.012095753375457\\
23	0.0114919120191121\\
24	0.0114915874902618\\
25	0.0114791822764769\\
26	0.0114791164300901\\
27	0.00942848486431508\\
28	0.00942818463158277\\
29	0.00753288540033795\\
30	0.00753278297038982\\
31	0.00752777202203293\\
32	0.0075272738516459\\
33	0.00720851011526472\\
34	0.00720817840874839\\
35	0.00712811894643316\\
36	0.00712798767233169\\
37	0.00598874820673185\\
38	0.00598869108086194\\
39	0.0051289427089822\\
40	0.00512891965034824\\
41	0.00503424481958879\\
42	0.00503424452216082\\
43	0.00494515674123981\\
44	0.00494492780082036\\
45	0.00487137514777922\\
46	0.00487037432500624\\
47	0.00421444940822556\\
48	0.00421409967359719\\
49	0.00395873561982647\\
50	0.00395846526079949\\
51	0.00386362980513069\\
52	0.00386353684236286\\
53	0.00377082854130412\\
54	0.00377030557373398\\
55	0.00372967042540982\\
56	0.00372896245992916\\
57	0.00351333695898237\\
58	0.00351271540544044\\
59	0.00330378411420591\\
};
\addlegendentry{$\frac{\Vert q - q_n\Vert}{\Vert q \Vert}$}

\addplot [color=mycolor2, dashed]
  table[row sep=crcr]{%
0	1\\
1	0.359805938935254\\
2	0.359781909976563\\
3	0.0545753225224157\\
4	0.0500049248468868\\
5	0.0171580345029261\\
6	0.017158024593463\\
7	0.00779911962046153\\
8	0.00778353165319985\\
9	0.0042627190538161\\
10	0.00426251172388405\\
11	0.00262415547666798\\
12	0.00262415415631536\\
13	0.00175172931821419\\
14	0.00175156803775569\\
15	0.00123912298873581\\
16	0.00123892084453071\\
17	0.000915256610565958\\
18	0.00091516079845879\\
19	0.000699140233566558\\
20	0.000699045948601238\\
21	0.000548607680110321\\
22	0.000548416146416228\\
23	0.000440031206785236\\
24	0.00043977995696157\\
25	0.000359558164702627\\
26	0.000359091518308983\\
27	0.000298257186097594\\
28	0.000297968409405683\\
29	0.000250882821298186\\
30	0.000250481873870689\\
31	0.000213482323544853\\
32	0.000213133469577126\\
33	0.000183854730360781\\
34	0.000183005373762078\\
35	0.000159403599160049\\
36	0.000159116464151412\\
37	0.000139579794036746\\
38	0.00013945920984066\\
39	0.000123354637878814\\
40	0.000123184854554119\\
41	0.000109956893484961\\
42	0.000109689916688902\\
43	9.88863166128095e-05\\
44	9.84504105634589e-05\\
45	8.96135256736368e-05\\
46	8.92384571188144e-05\\
47	8.18505811983025e-05\\
48	8.17956139497048e-05\\
49	7.59996218074654e-05\\
50	7.52597664068509e-05\\
51	7.07158604120756e-05\\
52	7.07106337064115e-05\\
53	6.70454953857494e-05\\
54	6.69344981723225e-05\\
55	6.44948263828488e-05\\
56	6.42298844051814e-05\\
57	6.28061113603602e-05\\
58	6.27910691174383e-05\\
59	6.22824346975996e-05\\
};
\addlegendentry{$ \frac{\Vert\phi-\phi_n\Vert}{\Vert \phi \Vert}$}

\addplot [color=black!30!mycolor2, mark size=0.8pt, mark=*, mark options={solid, fill=mycolor2, black!30!mycolor2}]
  table[row sep=crcr]{%
0	1\\
1	0.114746569346325\\
2	0.108715634236489\\
3	0.0450522534038761\\
4	0.0374498941113154\\
5	0.0181735389760626\\
6	0.0146773480144537\\
7	0.00978676663765104\\
8	0.00783460523447255\\
9	0.00502331305582069\\
10	0.00418318670877018\\
11	0.00215837737150043\\
12	0.00181254054959292\\
13	0.00119371870838012\\
14	0.00124273657563884\\
15	0.0013592880812847\\
16	0.00130839543498418\\
17	0.00126335206939511\\
18	0.00107346853540631\\
19	0.0008359911810919\\
20	0.000660255431805165\\
21	0.000369286463864993\\
22	0.000334743045209134\\
23	0.000307987996633321\\
24	0.000356998703881149\\
25	0.000414889000216105\\
26	0.000398806602529724\\
27	0.000385282073182484\\
28	0.00032799906449284\\
29	0.00025515368600473\\
30	0.00021397382639723\\
31	0.000137710698600638\\
32	0.000119900079294296\\
33	0.000148214466374314\\
34	0.000162383477899913\\
35	0.000174876144600143\\
36	0.000164410209801904\\
37	0.000162608212125621\\
38	0.000142256243929035\\
39	0.000111834274284489\\
40	9.22647874621909e-05\\
41	7.65407876379624e-05\\
42	7.59304873291416e-05\\
43	0.000100163528522023\\
44	0.000110763683328318\\
45	0.000120706775731265\\
46	0.000108558129344227\\
47	9.53491419425109e-05\\
48	7.46230170551248e-05\\
49	5.01243240728103e-05\\
50	4.66349602546242e-05\\
51	3.38326534582241e-05\\
52	5.85898716148313e-05\\
53	7.66496277901862e-05\\
54	8.57677161740406e-05\\
55	9.4280784157448e-05\\
56	8.54197595429984e-05\\
57	7.60067300788018e-05\\
58	5.93422600159246e-05\\
59	2.91370772041188e-05\\
};
\addlegendentry{$\frac{\Vert q - f(\phi_n)\Vert}{\Vert q \Vert}$}

\end{axis}

\begin{axis}[%
width=0.262\figurewidth,
height=\figureheight,
at={(0\figurewidth,0\figureheight)},
scale only axis,
xmin=0,
xmax=1,
xlabel style={font=\color{white!15!black}},
xlabel={$x$},
ymin=-0.2,
ymax=0.5,
axis background/.style={fill=white},
title style={font=\bfseries},
title={levelset function},
scaled ticks=false, xticklabel style={/pgf/number format/fixed},yticklabel style={/pgf/number format/fixed}
]
\addplot [color=mycolor1, dashed, forget plot]
  table[row sep=crcr]{%
0	0.847721129991393\\
0.01	0.812978988936933\\
0.02	0.778903514549138\\
0.03	0.745494706828011\\
0.04	0.71275256577355\\
0.05	0.680677091385756\\
0.06	0.649268283664628\\
0.07	0.618526142610167\\
0.08	0.588450668222373\\
0.09	0.559041860501246\\
0.1	0.530299719446785\\
0.11	0.50222424505899\\
0.12	0.474815437337863\\
0.13	0.448073296283402\\
0.14	0.421997821895608\\
0.15	0.39658901417448\\
0.16	0.371846873120019\\
0.17	0.347771398732225\\
0.18	0.324362591011097\\
0.19	0.301620449956637\\
0.2	0.279544975568842\\
0.21	0.258136167847715\\
0.22	0.237394026793254\\
0.23	0.21731855240546\\
0.24	0.197909744684332\\
0.25	0.179167603629871\\
0.26	0.161092129242077\\
0.27	0.143683321520949\\
0.28	0.126941180466489\\
0.29	0.110865706078694\\
0.3	0.0954568983575668\\
0.31	0.0807147573031059\\
0.32	0.0666392829153116\\
0.33	0.0532304751941841\\
0.34	0.0404883341397232\\
0.35	0.028412859751929\\
0.36	0.0170040520308014\\
0.37	0.00626191097634053\\
0.38	-0.00381356341145371\\
0.39	-0.0132223711325813\\
0.4	-0.0219645121870422\\
0.41	-0.0300399865748363\\
0.42	-0.0374487942959639\\
0.43	-0.0441909353504248\\
0.44	-0.050266409738219\\
0.45	-0.0556752174593466\\
0.46	-0.0604173585138075\\
0.47	-0.0644928329016017\\
0.48	-0.0679016406227292\\
0.49	-0.0706437816771901\\
0.5	-0.0727192560649844\\
0.51	-0.0741280637861119\\
0.52	-0.0748702048405728\\
0.53	-0.074945679228367\\
0.54	-0.0743544869494946\\
0.55	-0.0730966280039555\\
0.56	-0.0711721023917497\\
0.57	-0.0685809101128773\\
0.58	-0.0653230511673382\\
0.59	-0.0613985255551324\\
0.6	-0.0568073332762599\\
0.61	-0.0515494743307208\\
0.62	-0.045624948718515\\
0.63	-0.0390337564396426\\
0.64	-0.0317758974941035\\
0.65	-0.0238513718818977\\
0.66	-0.0152601796030252\\
0.67	-0.00600232065748609\\
0.68	0.0039222049547197\\
0.69	0.014513397233592\\
0.7	0.0257712561791312\\
0.71	0.0376957817913369\\
0.72	0.0502869740702094\\
0.73	0.0635448330157485\\
0.74	0.0774693586279543\\
0.75	0.0920605509068268\\
0.76	0.107318409852366\\
0.77	0.123242935464572\\
0.78	0.139834127743444\\
0.79	0.157091986688983\\
0.8	0.175016512301189\\
0.81	0.193607704580062\\
0.82	0.2128655635256\\
0.83	0.232790089137806\\
0.84	0.253381281416679\\
0.85	0.274639140362218\\
0.86	0.296563665974424\\
0.87	0.319154858253296\\
0.88	0.342412717198835\\
0.89	0.366337242811041\\
0.9	0.390928435089914\\
0.91	0.416186294035453\\
0.92	0.442110819647659\\
0.93	0.468702011926531\\
0.94	0.49595987087207\\
0.95	0.523884396484276\\
0.96	0.552475588763148\\
0.97	0.581733447708687\\
0.98	0.611657973320893\\
0.99	0.642249165599766\\
1	0.673507024545305\\
};
\addplot [color=mycolor2, forget plot]
  table[row sep=crcr]{%
0	0.376133385176628\\
0.01	0.366133409770264\\
0.02	0.356133435335724\\
0.03	0.34613346193177\\
0.04	0.336133489622002\\
0.05	0.326133518475361\\
0.06	0.316133548566705\\
0.07	0.306133579977456\\
0.08	0.296133612796337\\
0.09	0.286133647120208\\
0.1	0.276133683055019\\
0.11	0.266133720716907\\
0.12	0.256133760233447\\
0.13	0.246133801745093\\
0.14	0.236133845406851\\
0.15	0.226133891390212\\
0.16	0.2161339398854\\
0.17	0.206133991104005\\
0.18	0.196134045282069\\
0.19	0.186134102683722\\
0.2	0.17613416360549\\
0.21	0.166134228381418\\
0.22	0.156134297389195\\
0.23	0.146134371057528\\
0.24	0.136134449875054\\
0.25	0.126134534401189\\
0.26	0.116134625279431\\
0.27	0.106134723253769\\
0.28	0.0961348291890982\\
0.29	0.08613494409682\\
0.3	0.0761350691672347\\
0.31	0.066135205810927\\
0.32	0.0561353557121958\\
0.33	0.0461355208988296\\
0.34	0.0361357038343717\\
0.35	0.0261359075418147\\
0.36	0.0161361357719639\\
0.37	0.00613639323649742\\
0.38	-0.00386331406329249\\
0.39	-0.0138629783628802\\
0.4	-0.0238625894356185\\
0.41	-0.0338621335328176\\
0.42	-0.0438615917247162\\
0.43	-0.0538609372061878\\
0.44	-0.0638601307260436\\
0.45	-0.0738591124152337\\
0.46	-0.0838577862038338\\
0.47	-0.093855987589775\\
0.48	-0.103853409517265\\
0.49	-0.113849405304436\\
0.5	-0.123842339926807\\
0.51	-0.133826538275112\\
0.52	-0.143759924052693\\
0.53	-0.14596314088803\\
0.54	-0.136084040990588\\
0.55	-0.126104151012919\\
0.56	-0.116112402232157\\
0.57	-0.106116894297044\\
0.58	-0.0961197192861938\\
0.59	-0.0861216599036865\\
0.6	-0.0761230752052982\\
0.61	-0.0661241530547455\\
0.62	-0.056125001281249\\
0.63	-0.0461256861966962\\
0.64	-0.0361262508223736\\
0.65	-0.0261267242887688\\
0.66	-0.0161271270230879\\
0.67	-0.00612747377364165\\
0.68	0.0038722245449643\\
0.69	0.0138719596821787\\
0.7	0.023871725285451\\
0.71	0.0338715163841346\\
0.72	0.0438713290331059\\
0.73	0.0538711600612616\\
0.74	0.0638710068905701\\
0.75	0.0738708674036174\\
0.76	0.0838707398451328\\
0.77	0.093870622747742\\
0.78	0.103870514875267\\
0.79	0.113870415178919\\
0.8	0.123870322763086\\
0.81	0.133870236858355\\
0.82	0.14387015680004\\
0.83	0.153870082010951\\
0.84	0.163870011987453\\
0.85	0.173869946288113\\
0.86	0.183869884524386\\
0.87	0.193869826352922\\
0.88	0.203869771469181\\
0.89	0.213869719602102\\
0.9	0.223869670509623\\
0.91	0.233869623974906\\
0.92	0.243869579803134\\
0.93	0.253869537818785\\
0.94	0.263869497863301\\
0.95	0.273869459793086\\
0.96	0.283869423477784\\
0.97	0.293869388798784\\
0.98	0.303869355647931\\
0.99	0.313869323926394\\
1	0.323869293543685\\
};
\end{axis}
\end{tikzpicture}%

%% file: images/tikz/L2error1.tex
%
\definecolor{mycolor1}{rgb}{0.00000,0.44700,0.74100}%
\definecolor{mycolor2}{rgb}{0.85000,0.32500,0.09800}%
\definecolor{mycolor3}{rgb}{0.92900,0.69400,0.12500}%
\begin{tikzpicture}

\begin{axis}[%
width=0.98\figurewidth,
height=\figureheight,
at={(0\figurewidth,0\figureheight)},
scale only axis,
xmin=0,
xmax=35,
xlabel style={font=\color{white!15!black}},
xlabel={$n$ modes},
ymode=log,
ymin=1e-05,
ymax=1,
yminorticks=true,
ylabel style={font=\color{white!15!black}},
ylabel={relative error},
axis background/.style={fill=white},
legend style={legend cell align=left, align=left, fill=none, draw=none},
scaled ticks=false, xticklabel style={/pgf/number format/fixed},yticklabel style={/pgf/number format/fixed}
]
\addplot [color=mycolor1, mark=+, mark options={solid, mycolor1}]
  table[row sep=crcr]{%
0	1\\
1	0.208229443543112\\
2	0.102534111517672\\
3	0.0684607934553688\\
4	0.0525537994592058\\
5	0.0390562664701335\\
6	0.0314835452029856\\
7	0.0269628292186073\\
8	0.0219475880704862\\
9	0.0188429414737321\\
10	0.0163972664662912\\
11	0.013965614250307\\
12	0.0125463122580134\\
13	0.0109721803712042\\
14	0.00955260974393814\\
15	0.00883978008328678\\
16	0.00751239333211285\\
17	0.00653701996177867\\
18	0.00619338004598776\\
19	0.00522964408485424\\
20	0.00480151486331356\\
21	0.00442394922096749\\
22	0.00381990100342253\\
23	0.00353175709993112\\
24	0.00322741659242932\\
25	0.00284387635092358\\
26	0.00259873041169349\\
27	0.00242590484263579\\
28	0.00222875988078201\\
29	0.00221177107765771\\
30	0.0018469066254779\\
31	0.00169378260522157\\
32	0.00118759149876739\\
33	0.00106989788831319\\
34	0.000736697540171317\\
35	0.000718134579281819\\
36	0.000686160746658145\\
37	0.000402244487173694\\
38	0.00035469526584975\\
};
\addlegendentry{$\Vert q-q_n\Vert_2/\Vert q\Vert_2$}

\addplot [color=mycolor2, mark=x, mark options={solid, mycolor2}]
  table[row sep=crcr]{%
0	1\\
1	0.265510800527437\\
2	0.0463614181640376\\
3	0.0283715884758292\\
4	0.0199029658072694\\
5	0.0126201418660051\\
6	0.0118068407836256\\
7	0.00922624224162863\\
8	0.00722770731927204\\
9	0.00580362197791494\\
10	0.00468439293661942\\
11	0.00452342661533994\\
12	0.00362688541669926\\
13	0.00357160680174792\\
14	0.00352620886483618\\
15	0.003414557436661\\
16	0.00336364279101439\\
17	0.00333974173521294\\
18	0.00325883933378627\\
19	0.00323821649124755\\
20	0.00323636312864429\\
21	0.00323328260047271\\
22	0.00323555684310697\\
23	0.00323525501204207\\
24	0.00323540938117981\\
25	0.00323492470907774\\
26	0.0032351082490826\\
27	0.00323499042213585\\
28	0.00323481579946083\\
29	0.00323483255312862\\
30	0.00323472085101366\\
31	0.00323470372267496\\
32	0.0032347178590847\\
33	0.00323472310833984\\
34	0.00323470720702309\\
35	0.0032347257942111\\
36	0.00323474819740701\\
37	0.00323475213020157\\
38	0.0032347349939891\\
39	0.00323474067829452\\
};
\addlegendentry{$\Vert q-f(\phi_n)\Vert_2/\Vert q\Vert_2$}

\addplot [color=mycolor3, mark=o, mark options={solid, mycolor3}]
  table[row sep=crcr]{%
0	1\\
1	0.204567953291119\\
2	0.0443917947230602\\
3	0.0232309619741152\\
4	0.0085709403775077\\
5	0.00731862287415803\\
6	0.00489921557644142\\
7	0.00283262683278337\\
8	0.00178617558199724\\
9	0.00164870807510567\\
10	0.00120265155748921\\
11	0.00111233314508517\\
12	0.000823138940726148\\
13	0.000677207015357723\\
14	0.000553609434910426\\
15	0.000528998406168389\\
16	0.00045021550101092\\
17	0.000379150226217814\\
18	0.000338482402789455\\
19	0.000280990340105079\\
20	0.000254682920176477\\
21	0.000242397249647828\\
22	0.000218394872009236\\
23	0.000207289773748153\\
24	0.00016621041992349\\
25	0.000154863734452438\\
26	0.000120363422383624\\
27	0.000109332163294977\\
28	8.97499587945491e-05\\
29	8.94511649719974e-05\\
30	8.04497242724231e-05\\
31	7.50541169400655e-05\\
32	7.2479741508773e-05\\
33	6.90021502406033e-05\\
34	6.22861791475555e-05\\
35	6.12830318216549e-05\\
36	5.66541279019915e-05\\
37	5.29509644899316e-05\\
38	4.88074401371284e-05\\
};
\addlegendentry{$\Vert\phi-\phi_n\Vert_2/\Vert\phi\Vert_2$}

\addplot [color=black, dashed, line width=2.0pt, forget plot]
  table[row sep=crcr]{%
1	0.00323473166550014\\
39	0.00323473166550014\\
};
\node[right, align=left]
at (axis cs:2,0.002) {$\Delta f$};
\end{axis}
\end{tikzpicture}%